\documentclass[11pt]{siamart220329}
\usepackage{times,amsmath}%
\usepackage{comment}

\usepackage{bm}
\usepackage{amssymb,amsfonts,amscd,cleveref}
\usepackage{graphicx,latexsym,color}
\usepackage{hyperref}
\hypersetup{
    colorlinks=true,
    linkcolor=blue,
    filecolor=magenta,      
    urlcolor=cyan,
    pdftitle={Overleaf Example},
    pdfpagemode=FullScreen,
    }

\usepackage{multirow}
\usepackage{wrapfig}
\usepackage{soul}
\DeclareUnicodeCharacter{0301}{\'{}}

\usepackage{subfig}
\newcommand{\sign}{\text{sign}}
\newtheorem{remark}{\textbf{Remark}}[section]
\def \bes{\begin{eqnarray}}
\def \ees{\end{eqnarray}}
\def \bns{\begin{eqnarray*}}
\def \ens{\end{eqnarray*}}

\newenvironment{eqa}{\begin{equation}%
  \begin{array}{rcl}}{\end{array}\end{equation}}
\newcommand\beqa{\begin{eqa}}
\newcommand\eeqa{\end{eqa}}
\usepackage{tikz}
\usetikzlibrary{arrows}

\usepackage{pgfplots}
\pgfplotsset{compat=1.17}
\usepgfplotslibrary{external}
\usepackage[bottom=1.2in,top=1in,left=1in,right=1in]{geometry}

\title{Stability and Robustness of Time-discretization Schemes for the Allen-Cahn Equation via Bifurcation and Perturbation Analysis}
\author{%
  Wenrui Hao\thanks{Department of Mathematics, The Pennsylvania State University, State College, PA, 16802, USA(\email{wxh64@psu.edu}).}
\and
  Sun Lee\thanks{Department of Mathematics, The Pennsylvania State University, State College, PA, 16802, USA(\email{skl5876@psu.edu}).}
    \and
Xiaofeng Xu\thanks{Applied Mathematics and Computational Sciences, King Abdullah University of Science and Technology, 23955, Saudi Arabia(\email{xiaofeng.xu@kaust.edu.sa}).}
    \and
  Zhiliang Xu\thanks{Applied and Computational Mathematics and Statistics Department, University of Notre Dame, Notre Dame, IN 46556, USA(\email{zxu2@nd.edu
}).}}

\begin{document}
\maketitle
\begin{abstract}
The Allen-Cahn equation is a fundamental model for phase transitions, offering critical insights into the dynamics of interface evolution in various physical systems. This paper investigates the stability and robustness of frequently utilized time-discretization numerical schemes for solving the Allen-Cahn equation, with focuses on the Backward Euler, Crank-Nicolson (CN), convex splitting of modified CN, and Diagonally Implicit Runge-Kutta (DIRK) methods. Our stability analysis reveals that the Convex Splitting of the Modified CN scheme exhibits unconditional stability, allowing greater flexibility in time step selection, while the other schemes are conditionally stable. Additionally, our robustness analysis highlights that the Backward Euler method converges to correct physical solutions regardless of initial conditions. In contrast, the other methods studied in this work show sensitivity to initial conditions and may converge to incorrect physical solutions if the initial conditions are not carefully chosen. This study introduces a comprehensive approach to assessing stability and robustness in numerical methods for solving the Allen-Cahn equation, providing a new perspective for evaluating numerical techniques for general nonlinear differential equations.
\end{abstract}

\section{Introduction}
The Allen-Cahn equation, a fundamental partial differential equation (PDE) in the field of phase transitions, describes the process of phase separation in multi-component alloy systems. Its significance extends to numerous applications in materials science \cite{heida2012development}, image processing \cite{benevs2004geometrical, lee2019image}, and other areas requiring the modeling of interface dynamics \cite{aihara2019multi, grant1983temperature}. Due to the equation's nonlinearity and the presence of diffuse interface in solutions, developing robust and stable numerical schemes is a long-lasting challenge for accurate simulations \cite{du2018stabilized, feng2003numerical, jeong2016comparison, shen2010numerical}.  As important as spatial discretization, time discretization is crucial since it also directly determines the efficiency and accuracy of the numerical schemes \cite{du2019maximum, du2020time, li2010unconditionally,  liu2003phase, zhang2009numerical}. Below is an incomplete list of frequently utilized explicit and implicit time discretization methods for solving the Allen-Cahn equation and other  phase field models:

\begin{itemize}
    \item \textbf{Explicit Methods}
    \begin{itemize}
        \item \textit{Forward Euler Method}: This first-order method approximates the time derivative using a simple forward difference. It is conditionally stable and often requires very small time steps, especially for stiff problems like the Allen-Cahn equation \cite{feng2003numerical, jeong2016comparison}.
        \item \textit{Runge-Kutta Methods}: Higher-order explicit methods, such as the fourth-order Runge-Kutta, can be used to improve accuracy while still being conditionally stable \cite{zhang2021numerical}. These methods are rarely used due to the stringent time-step restrictions imposed by stability considerations.
    \end{itemize}

    \item \textbf{Implicit Methods}
    \begin{itemize}
        \item \textit{Backward Euler Method}: A first-order implicit method that is unconditionally stable and well-suited for stiff problems \cite{feng2003numerical, jeong2016comparison}. However, it requires solving a nonlinear system at each time step.
        \item \textit{Crank-Nicolson (CN) Method}: This second-order implicit method is based on the trapezoidal rule and offers a good balance between accuracy and stability \cite{feng2013stabilized, hou2017numerical, zhang2009numerical}. It is also unconditionally stable but, like the backward Euler method, requires solving a nonlinear system at each step.
   \item \textit{Diagonally Implicit Runge-Kutta (DIRK) Methods}: These methods are a subclass of implicit Runge-Kutta methods where the coefficient matrix is lower triangular with equal diagonal elements \cite{zhang2021preserving}. This structure simplifies the implementation by allowing a step-by-step solution of implicit equations, which improves stability and accuracy while maintaining reasonable computational costs.

    \end{itemize}

    \item \textbf{Semi-Implicit Methods}
    \begin{itemize}
        \item \textit{Semi-Implicit Spectral Deferred Correction (SISDC) Method}: This method iteratively corrects the solution using both implicit and explicit updates, improving stability and accuracy \cite{liu2015stabilized}.
        \item \textit{Semi-Implicit Backward Euler Method}: This approach involves treating the stiff linear terms implicitly while handling the nonlinear terms explicitly, reducing the complexity of solving fully implicit equations. \cite{chen1998applications,shen2010numerical,spigler1995convergence}
        \item {\textit{Convex-splitting Method}: This method explores careful splitting of the nonconvex term into the difference of two convex terms, making them implicit and explicit, respectively. See \cite{guillen2013linear,eyre1998unconditionally, graser2013time, guillen2014second}.} 
    \end{itemize}
\end{itemize}

The choice of time discretization method for the Allen-Cahn equation and other phase field models requires careful consideration of stability, accuracy, and computational efficiency \cite{xu2019stability}. 
Implicit and semi-implicit methods are often favored for their desired stability properties, and are particularly suited for stiff problems. Notably, these methods require solving nonlinear systems.

 Fourier or
energy methods are often used to analyze stability conditions for linear schemes of linear partial differential equations with constant coefficients. Yet few have become known about the stability of numerical schemes
for nonlinear equations. The study in \cite{xu2023lack} indicates that many numerical schemes, except for the backward Euler method, may experience convergence issues unless the time step size is exceedingly small. Motivated by this work, we introduce in this paper the following concepts of stability and robustness for numerical schemes designed to solve the Allen-Cahn equation.
\begin{definition}
     Stability is defined as the uniqueness of $\phi^{n+1}$ given $\phi^n$, revealing the upper bound for the step size of the numerical scheme which is the stability condition;
     \end{definition}
\begin{definition}     
     Robustness is defined as the uniqueness of $\phi^n$ given $\phi^{n+1}$, indicating the numerical scheme's accuracy in converging to the physical solution.
\end{definition}

While the stability of numerical schemes is a center of concern in analyzing them, the importance of robustness of the nonlinear schemes is sometimes overlooked. In fact, the robustness here indicates the sensitivity of the (nonlinear) schemes to the initial guess used to solve them in each time-stepping. Thus a numerical solution computed using a scheme suffering from robustness issues may converge to a wrong solution. 

 Both stability and robustness require the application of bifurcation theory \cite{chow2012methods, iooss2012elementary, kuznetsov1998elements, sattinger2006topics} and perturbation analysis \cite{bonnans2013perturbation, ho2012perturbation}, powerful mathematical tools that examine the behavior of solution structures. Bifurcation analysis allows one to explore the uniqueness of the numerical solutions and identify critical points where qualitative changes in the solution structure occur. Perturbation analysis provides insight into how small perturbations affect the structure of trivial solutions. Together, these analyses form a rigorous framework to evaluate the performance of different numerical schemes.

In this paper, we examine both the stability and robustness of several time-discretization numerical schemes that are commonly used for the Allen-Cahn equation. They include the Backward Euler method, Crank-Nicolson method, and Runge-Kutta methods. Our findings reveal both essential stability conditions and sensitivity to initial conditions based on robustness, guiding the development of efficient and reliable computational methods for simulating the Allen-Cahn equation. Through this investigation, we aim to develop a general framework based on stability and robustness in the numerical treatment of phase field models.

The rest of the paper is organized as follows: In Section 2, we present the Allen-Cahn equation considered in this study. Section 3 explores the stability and robustness of the Backward Euler scheme. In Section 4, we extend this analysis to the Crank-Nicolson scheme. Section 5 examines the convex splitting of the modified Crank-Nicolson scheme. In Section 6, we study the Diagonally Implicit Runge-Kutta method. Finally, we conclude our work in Section 7.

\section{Problem Setup}
We consider the following time-dependent  Allen-Cahn equation on domain $[-1,1]^d=\Omega \subset \mathbb{R}^d$:
\begin{equation}   
\begin{cases}
\phi_t-\triangle\phi+\frac{1}{\varepsilon^2}(\phi^3-\phi)=0\hbox{~,~~~~~~} x\in \Omega~, ~ t \in [0, T]~, \\ 
\phi(x, 0) = \phi_0(x) \hbox{~,~~~~~~} x\in \Omega~, \\ 
 \frac{\partial \phi(x,t)}{\partial \mathbf{n}}=0~\hbox{~,~~~~~~} x\in \partial\Omega~, ~ t \in [0, T]~
\end{cases}
\label{AC}
\end{equation}
where $\varepsilon$ is a small positive constant representing the thickness of the diffuse interface,  $\phi_0(x)$ is given and $\mathbf{n}$ is the unit outward normal vector to $\partial \Omega$, $x=(x_1,x_2,\cdots,x_d)^T$. {It is well known that the Allen-Cahn equation possesses two stable steady state solutions $\phi(x) = \pm 1$, respectively. }

Our objective is to investigate numerical methods for solving the Allen-Cahn equation. We consider nonlinear schemes given in the general formula $ F(\phi^n(x), \phi^{n+1}(x)) = 0$, where $ \phi^n(x)$ represents the solution at time step $n$, and $\phi^{n+1}(x)$ represents the solution at the next time step $n+1$.

Specifically, we explore the following two distinct aspects:
\begin{itemize}
    \item[1)] {\bf Stability:} For a given $\phi^n(x)$, we analyze the uniqueness of $\phi^{n+1}(x)$ through bifurcation analysis, establishing a stability condition for the numerical scheme.
\item[2)] {\bf Robustness:} For a given $\phi^{n+1}(x)$, we investigate the solution landscape of $ \phi^n(x)$ to assess the robustness of the scheme to the initial guess.
\end{itemize}

\section{Backward Euler Scheme}\label{BEuler}
The backward Euler scheme for the Allen-Cahn Eq.~\eqref{AC} reads as:
\begin{equation}\label{fis1st}
\begin{cases}
\frac{\phi^{n+1}-\phi^n}{\Delta t}-\triangle\phi^{n+1}+\frac{1}{\varepsilon^2}\big((\phi^{n+1})^3-\phi^{n+1}\big)=0~,
    \\   
 \frac{\partial \phi^n(x)}{\partial \mathbf{n}}= \frac{\partial \phi^{n+1}(x)}{\partial \mathbf{n}}=0~\hbox{~,~~~~~~} x\in \partial\Omega~,
\end{cases}
\end{equation}
where $\Delta t$ is the time step size. We will derive its stability condition through bifurcation analysis and its robustness through perturbation analysis in the following subsection.



\subsection{The stability condition via bifurcation analysis}
In this section, we begin with the constant solution of Eq.~\eqref{fis1st} and then study the bifurcation analysis of the constant solution with respect to the parameter $\varepsilon$.

First, we examine a constant solution case of Eq.~\eqref{fis1st}, where $\phi^n \equiv r$ and $\phi^{n+1} \equiv c$, satisfying the equation:
\[
\frac{c-r}{\Delta t} + \frac{1}{\varepsilon^2}(c^3 - c) = 0~.
\]

Next, we consider a general perturbed case by introducing $\phi^n(x) = r + \delta f(x)$ and $\phi^{n+1}(x) = c + \delta \psi(x)$, with $\delta \in \mathbb{R}$, $|\delta| \ll 1$, and $f$ and $\psi$ being some smooth functions. 
Substituting these functions into Eq.~\eqref{fis1st} and retaining the linear term in $\delta$, we obtain
\begin{equation}
\frac{\delta \psi -\delta f}{\Delta t}-\delta \triangle\psi+\frac{1}{\varepsilon^2}\left(3c^2- 1 \right)\delta \psi+O(\delta^2)=0~.
\end{equation}

Rearranging this equation, we obtain
\begin{equation}\label{BE}
-\triangle\psi + \left
(\frac{ 1 }{\Delta t}+ \frac{3c^2- 1}{\varepsilon^2} \right) \psi  + O(\delta) = \frac{f}{\Delta t}~.
\end{equation}


We can assess the uniqueness of $\phi^{n+1}$ by studying the solution structure of Eq.~\eqref{BE} after dropping the $O(\delta)$ term.
In this case, $f(x)$ is a given function, and we focus on examining the homogeneous part of Eq.~\eqref{BE}, namely,
\begin{equation}
-\triangle\psi + \left( \frac{ 1 }{\Delta t}+ \frac{3c^2- 1}{\varepsilon^2} \right) \psi = 0, \label{Linear1}
\end{equation} 
which is a Helmholtz equation. We arrive at the following proposition.

\begin{proposition}\label{prop1}
Bifurcations of $\psi$ in Eq.~\eqref{Linear1} occur when $1 - 3c^2 > 0$ and $\Delta t > \frac{\epsilon^2}{1 - 3 c^2} \geq \epsilon^2 $. The bifurcation points and corresponding eigenfunctions are as follows:
\begin{align}
\epsilon^2& =\frac{1 - 3c^2}{\frac{1}{\Delta t}+\sum\limits_{i=1}^d\pi^2k_i^2}, && \psi(x)=\prod_{i=1}^d A_i(x_i), &&A_i(x_i)=\begin{cases}
\cos(\pi k_ix_i),~ \quad k_i=0,1,2,\cdots,
    \\   
\sin(\pi k_ix_i),~ \quad k_i=\frac{1}{2},\frac{3}{2},\frac{5}{2},\cdots,
\end{cases}
\end{align}
where the function $A_i(x_i)$ can either be $\cos(\pi k_ix_i)$ or $\sin(\pi k_ix_i)$ depending the values of $k_i$.
\end{proposition}

\begin{proof}
By the method of separation of variables, we let $\psi(x)=\prod_{i=1}^d A_i(x_i)$. Substituting this into Eq.~\eqref{Linear1} yields:
\begin{align*} 
\left(\frac{1}{\Delta t}+\frac{3c^2 - 1}{\epsilon^2}\right)\prod_{i=1}^d A_i(x_i) - \sum_{i=1}^d \left( \partial_{x_i}^2A_i(x_i) \prod_{j \neq i} A_j(x_j) \right) = 0.
\end{align*}
For the trivial solution, this equation holds for any constant $\epsilon$. 
For the non-trivial case where $A_i \neq 0$, we divide both sides by $\prod_{i=1}^d A_i(x_i)$ to obtain:
\begin{align}\label{eq:prop-1} 
\sum_{i=1}^d \frac{\partial_{x_i}^2A_i(x_i)}{A_i(x_i)} = \left(\frac{1}{\Delta t}+\frac{3c^2 - 1}{\epsilon^2}\right).
\end{align}
Each term $\frac{\partial_{x_i}^2A_i(x_i)}{A_i(x_i)}$ in Eq.~\eqref{eq:prop-1} is constant because the sum $\sum_{i=1}^d \frac{\partial_{x_i}^2A_i(x_i)}{A_i(x_i)}$ remains constant regardless of changes in the variables $x_i$. 
This reduces our problem to a 1D eigenvalue problem with the Neumann boundary condition:
\begin{align}\label{1dEGP}
    \partial_{x_i}^2A_i(x_i)=c_iA_i(x_i),&& \left.\partial_{x_i}A_i(x_i) \right|_{x_i=-1}=\left.\partial_{x_i}A_i(x_i) \right|_{x_i=1}=0,&& i=1,2,\cdots,d,
\end{align}
where $\sum_{i=1}^dc_i=\left(\frac{1}{\Delta t}+\frac{3c^2 - 1}{\epsilon^2}\right)$. 

The eigenvalue $c_i$ and corresponding eigenfunction $A_i(x_i)$ for \eqref{1dEGP} can be easily computed, and are given by: 
\begin{align}
c_i &= k_i^2\pi^2 , && A_i(x_i)=\begin{cases}
\cos(\pi k_ix_i),~ \quad k_i=0,1,2,\cdots,
    \\   
\sin(\pi k_ix_i),~ \quad k_i=\frac{1}{2},\frac{3}{2},\frac{5}{2},\cdots,
\end{cases}
\end{align}
where the eigenfunction $A_i(x_i)$ can either be $\cos(\pi k_ix_i)$ or $\sin(\pi k_ix_i)$ depending the values of $k_i$.


To satisfy Eq.~\eqref{Linear1}, we also require $\sum_{i=1}^d c_i = -\sum_{i=1}^d \pi^2 k_i^2 = \left(\frac{1}{\Delta t} + \frac{3c^2 - 1}{\epsilon^2}\right)$. 
Thus, we have the bifurcation point $\epsilon^2 = \frac{1 - 3c^2}{\frac{1}{\Delta t} + \sum_{i=1}^d \pi^2k_i^2 }$ corresponding to the eigenfunction $\prod_{i=1}^d A_i(x_i)$.
\end{proof}

From Proposition \ref{prop1}, we can see that, if $\Delta t \le \epsilon^2$, the solution to Eq.~\eqref{Linear1} will exclusively exhibit the trivial solution $\psi = 0$. Consequently, we can conclude that the particular solution for Eq.~\eqref{BE} is unique, and $\phi^{n+1}$ is also unique while satisfying the stability condition $\Delta t \le \epsilon^2$. Thus, the backward Euler scheme is stable when $\Delta t \le \epsilon^2.$

\begin{remark}
    For the robustness analysis of the backward Euler scheme, we examine the uniqueness of $\phi^n$ given $\phi^{n+1}$. As Eq.~\eqref{fis1st} is linear with respect to \textcolor{black}{ $\phi^{n}$}, establishing uniqueness is straightforward.
\end{remark}




\section{Crank-Nicolson Scheme}
The Crank-Nicolson scheme for the Allen-Cahn Eq.~\eqref{AC}  reads as 
\begin{equation} \label{cn scheme}
\begin{cases}
    \frac{\phi^{n+1} - \phi^n}{\Delta t} -\frac{1}{2}(\triangle\phi^{n+1} + \triangle\phi^{n}) + \frac{1}{2\epsilon^2}\big( (\phi^{n+1})^3 - \phi^{n+1}\big) + \frac{1}{2\epsilon^2}\big( (\phi^{n})^3 - \phi^{n}\big) = 0~,
    \\
\frac{\partial \phi^n(x)}{\partial \mathbf{n}}= \frac{\partial \phi^{n+1}(x)}{\partial \mathbf{n}}=0~\hbox{~,~~~~~~} x\in \partial\Omega~.
\end{cases}
\end{equation}

\subsection{The stability condition via bifurcation analysis}
In this section, we investigate the uniqueness of $\phi^{n+1}$ for any given $\phi^n$ and derive the associated stability condition. We start by considering a perturbed setup with respect to the trivial solution. 

Here we define $\phi^{n+1}=c + \delta \psi(x)$ and $\phi^n=r + \delta f(x)$, with $c$ and $r$ representing constant solutions of Eq.~\eqref{cn scheme}, and $f(x)$ being a given function with $|\delta| \ll 1$. Plugging these expressions into Eq.~\eqref{cn scheme}, we obtain:
\begin{equation}
\frac{c + \delta \psi - r - \delta f}{\Delta t} -\frac{1}{2}(\delta \triangle\psi +\delta \triangle f) + \frac{1}{2\epsilon^2}\big( (c + \delta \psi)^3 - c - \delta \psi\big) + \frac{1}{2\epsilon^2}\big( ( r + \delta f)^3 - r - \delta f\big) = 0~.
\end{equation}

This equation is simplified as
\begin{equation}\label{parti_equ1}
-\frac{1}{2}\triangle\psi+\left(\frac{1}{\Delta t} + \frac{3c^2 - 1}{2\epsilon^2}\right)\psi + O(\delta) = G(x)~,
\end{equation}
where $G(x)$ is given by:

\begin{equation} 
G(x) = \frac{f(x)}{\Delta t} + \frac{1}{2} \triangle f(x) - \frac{1}{2\epsilon^2}\big( 3r^2f(x) - f(x)\big)~.
\end{equation}
Considering the homogeneous case of Eq.~\eqref{parti_equ1} and dropping the $O(\delta)$ term, we have:
\begin{equation}\label{cnho}
-\frac{1}{2}\triangle\psi+\left(\frac{1}{\Delta t} + \frac{3c^2 - 1}{2\epsilon^2}\right)\psi  = 0~.
\end{equation}
Then we can deduce the following bifurcation results:

\begin{proposition}\label{prop2}
Bifurcations of $\psi$ in Eq.~\eqref{cnho} occur when $1 - 3c^2 > 0$ and $\Delta t > \frac{2\epsilon^2}{1 - 3 c^2} \geq 2\epsilon^2$. The bifurcation points and corresponding eigenfunctions are as follows:
\begin{align}
\epsilon^2& =\frac{1 - 3c^2}{\frac{2}{\Delta t}+\sum\limits_{i=1}^d\pi^2k_i^2}, && \psi(x)=\prod_{i=1}^d A_i(x_i), &&A_i(x_i)=\begin{cases}
\cos(\pi k_ix_i),~ \quad k_i=0,1,2,\cdots,
    \\   
\sin(\pi k_ix_i),~ \quad k_i=,\frac{3}{2},\frac{5}{2},\cdots,
\end{cases}
\end{align}
where the function $A_i(x_i)$ can either be $\cos(\pi k_ix_i)$ or $\sin(\pi k_ix_i)$ depending the values of $k_i$.
\end{proposition}

From Proposition \ref{prop2}, we can see that, if $\Delta t \le 2\epsilon^2$, the solution to Eq.~\eqref{cnho} will exclusively exhibit the trivial solution $\psi = 0$. Consequently, we can conclude that the particular solution for Eq.~\eqref{parti_equ1} is unique, and $\phi^{n+1}$ is also unique while satisfying the stability condition $\Delta t \le 2\epsilon^2$. The proof proceeds by employing similar computations as those used in Proposition \ref{prop1}.

\subsection{The robustness analysis}
In this section, we delve into the solution space of $\phi^n$ given a  $\phi^{n+1}$. Initially, we investigate the trivial solutions; subsequently, we apply perturbation analysis to examine these trivial solutions.

\subsubsection{Trivial solution analysis}\label{3.2.1}
We consider $\phi^{n+1}\equiv c$ and $\phi^{n}\equiv r$ and rewrite Eq.~\eqref{cn scheme} as:
\begin{equation}\label{cr}
    \frac{c -  r }{\Delta t} + \frac{1}{2\epsilon^2}\left( c^3 - c \right) + \frac{1}{2\epsilon^2}\left( r^3 -  r \right) = 0~.
\end{equation}
By fixing $c$, if $-27\left(\frac{-2\epsilon^2 c}{\Delta t}-c^3+c\right)^2+4\left(\frac{2\epsilon^2}{\Delta t}+1\right)^3>0$,  we have 3 different real solutions for $r$ by the discriminant of cubic polynomial. If it is equal to $0$, then we have multiple real solutions. If it is less than $0$, we have complex solutions and a real solution.

First, we let $c=0$ and solve Eq.~\eqref{cr} to get two non-zero roots for $r$ and denote them as $\pm r_1$, where $r_1:=\sqrt{1 + \frac{2 \epsilon^2}{\Delta t}}$. 

Then we have the following results:
\begin{enumerate}
    \item $c = 0$ if and only if $r = 0$, $\pm r_1$;
    \item  $c > 0$ if and only if $ r \in \left(0, r_1\right)\cup \left(-\infty, -r_1\right)$;
    \item $c < 0$ if and only if $ r \in \left(-r_1, 0\right) \cup \left(r_1, \infty\right)$.
    \end{enumerate}
Therefore, given a time step size $\Delta t \leq 2\epsilon^2$ satisfying the stability condition, to have $\text{sign}(c) = \text{sign}(r)$ hold, we must have $|r| \leq r_1$;  to have $\text{sign}(c) \neq \text{sign}(r)$ hold, we must have $|r| > r_1$.
The above analysis is discussed in Theorem 3.2 of \cite{xu2023lack} to show that the Crank-Nicolson method may converge to a wrong steady-state solution. We carry out this further to obtain a more insightful convergence pattern.

Next, we compute roots by solving Eq.~\eqref{cr} with $c=r_1$. Since \begin{equation}
    -27\left(\frac{-2\epsilon^2 r_1}{\Delta t}-r_1^3+r_1\right)^2+4\left(\frac{2\epsilon^2}{\Delta t}+1\right)^3<0~,
\end{equation} we can obtain only one negative real solution and denote it as $-r_2$. Since $c=r_1>0$, we have $-r_2\in(-\infty,-r_1)$ which implies $r_1<r_2$. Moreover, notice that the value of $-27\left(\frac{-2\epsilon^2 c}{\Delta t}-c^3+c\right)^2+4\left(\frac{2\epsilon^2}{\Delta t}+1\right)^3$  decreases if $c$ is getting larger than $r_1$. {Thus we can find a unique sequence of $r_n$ by repeating this process with $c=-r_{i-1}$.}  This gives us a pattern of signs of trivial solutions at consecutive time steps.  We can conclude that:
\begin{enumerate}
    \item  $\phi^{n+1} > 0$ if and only if $ \phi^n \in \left(0, r_1\right)\cup \left(-\infty, -r_1\right)$;
        \item $\phi^{n+2} > 0$ if and only if $ \phi^n \in (0, r_1) \cup (-r_2,-r_1) \cup (r_2, \infty) $;

    \item $\phi^{n+1} < 0$ if and only if $ \phi^n \in \left(-r_1, 0\right) \cup \left(r_1, \infty\right)$;
            \item $\phi^{n+2} < 0$ if and only if $ \phi^n \in (-r_1,0) \cup (r_1,r_2) \cup (-\infty,-r_2) $.

    \end{enumerate}

\begin{figure}[ht]
    \centering
  \includegraphics[width=6in]{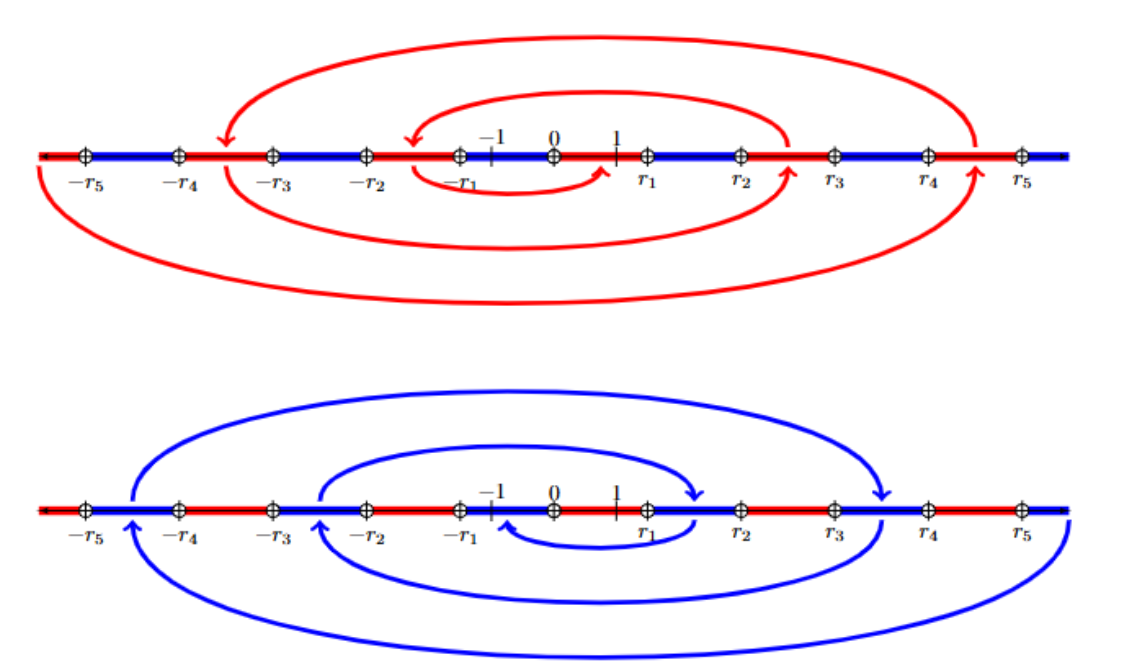}
    \caption{Visualizing Convergence Intervals of CN Scheme in Eq.~\eqref{cn scheme}. If initial conditions are chosen in red regions, CN Scheme eventually converges to $1$, while the initial conditions  chosen in  blue regions lead the CN Scheme to converge to $-1$. The values of $r_n$ with different $\frac{\Delta t}{2 \epsilon^2}$ are shown in Table~\ref{Table_CN}. }
    \label{fig:rc}
\end{figure}

Now we have generated a sequence of $r_i$ by solving Eq.~\eqref{cr} with setting $c=-r_{i-1}$. The numerical values of $r_i$ are presented in  Table~\ref{Table_CN} for different values of $\frac{\Delta t}{2\epsilon^2}$. The convergence intervals are shown in Fig.~\ref{fig:rc}. If one chooses the initial condition $\phi^n$ from the interval $[r_i,r_{i+1}]$, the CN scheme  converges to $(-1)^i\times [0,r_1]$ after $i$ time-stepping.

\begin{table}[ht]
\centering
\begin{tabular}{|c||c|c|c|c|}
\hline
$\frac{\Delta t}{2\epsilon^2}$ & $r_1$ & $r_2$ & $r_3$ & $r_4$\\
\hline
\hline
0.001 & 31.639 & 48.124 & 60.363 & 70.53 \\
\hline
0.01 & 10.05 & 15.256 & 19.123 & 22.335 \\
\hline
0.1 & 3.317 & 4.942 & 6.152 & 7.159 \\
\hline
0.25 & 2.236 & 3.243 & 3.996 & 4.625 \\ 
\hline
0.5 & 1.732 & 2.421 & 2.941 & 3.377 \\
\hline
\end{tabular}
\caption{Convergence Interval points $r_i$ for the CN Scheme in Eq.~\eqref{cn scheme} with different $\frac{\Delta t}{2\epsilon^2}$.}\label{Table_CN}
\end{table}

\subsubsection{Perturbation Analysis}\label{PertAna}

In this section, we delve into non-trivial solutions by perturbing the trivial solutions analyzed in the previous section. Specifically, we define $\phi^{n+1} = c + \delta f(x)$ and $\phi^{n} = r + \delta \psi(x)$. Our objective is to investigate $\phi^n$ for a given $\phi^{n+1}$. In this case, $f(x)$ is a given perturbation function, such as $\cos(k\pi x)$ {in the 1D case}. We need to solve for $\psi(x)$. To achieve this, we substitute these functions into Eq.~\eqref{cn scheme}, resulting in:
%
\begin{equation} 
    -\frac{1}{2} \triangle\psi +\left(-\frac{1}{\Delta t}+\frac{3r^2-1}{2\epsilon^2}\right)\psi +O(\delta)= G(x)~,
\end{equation}
where $G(x)$ is defined as
\begin{equation} 
    G(x)=-\frac{   f(x)}{\Delta t}+\frac{1}{2} \triangle f(x)-\frac{1}{2\epsilon^2}\big( 3c^2f(x) -  f(x)\big)~.
\end{equation}
To be more specific {for the 1D case}, when we choose $f(x) = \cos(k\pi x_1)$ for $k \in \mathbb{N}^+$, we find that $\psi = B \cos(k \pi x_1)$, with the coefficient $B$ given by 
\begin{equation}
    B = - \frac{\frac{3 c^2}{\epsilon^2} - \frac{1}{\epsilon^2} + \frac{2}{\Delta t} + (k \pi)^2}{ \frac{3 r^2}{\epsilon^2} - \frac{1}{\epsilon^2} - \frac{2}{\Delta t} + (k \pi)^2}~.
\end{equation}
For the 2D case, we choose $f(x) = \cos(k\pi x_1) \cos(l \pi x_2)$ for $k, l \in \mathbb{N}^+$, and have that $\psi = B \cos(k\pi x_1) \cos(l\pi x_2)$, with the coefficient $B$ given by 
\begin{equation}
    B = - \frac{\frac{3 c^2}{\epsilon^2} - \frac{1}{\epsilon^2} + \frac{2}{\Delta t} + (k \pi)^2+(l \pi)^2}{ \frac{3 r^2}{\epsilon^2} - \frac{1}{\epsilon^2} - \frac{2}{\Delta t} + (k \pi)^2+(l \pi)^2}~.
\end{equation}

\begin{figure}[!ht]
\centering
  \includegraphics[width=2in]{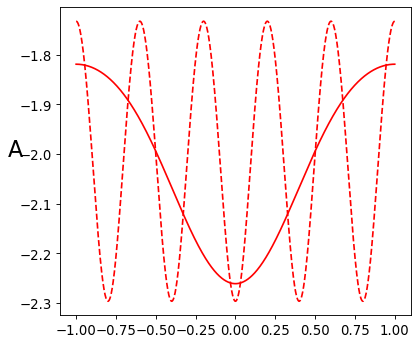}
  \includegraphics[width=2in]{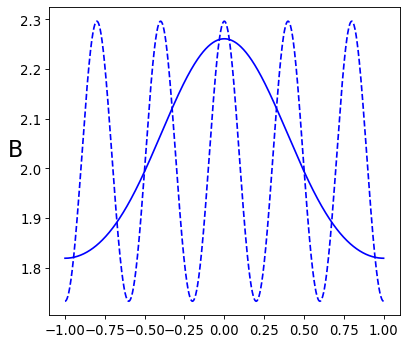}
  \includegraphics[width=2.2in]{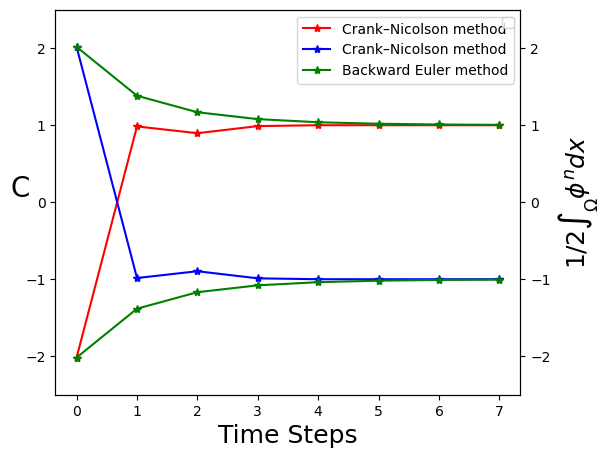}
  \caption{The solutions of $\phi^{n} \approx r + \delta B\cos(k\pi x_1)$ of the CN scheme for $\phi^{n+1}=c+\delta \cos(k\pi x_1)$ with $|\delta|=0.5$ are depicted in A and B panels. Here the solid curve is for $k=1$, the dashed curve is for $k=5$. The parameters are chosen as in A $(c=0.984375)$, ($r=-1.99310$) and B $(c=-0.984375)$, ($r=1.99310$), $\epsilon=0.1$, and $\Delta t=0.01$. These $r$ values are chosen from Table~\ref{Table_CN} in interval $(r_{1},r_{2})$. In panel C, the CN scheme jumps to a different solution with the initial conditions of $\phi^n$, and ultimately converges to an incorrect solution. Conversely, the backward Euler scheme converges to a solution near $\phi^n$.}\label{compare}
\end{figure}

\begin{figure}[!ht]
\centering
  \includegraphics[width=4in]{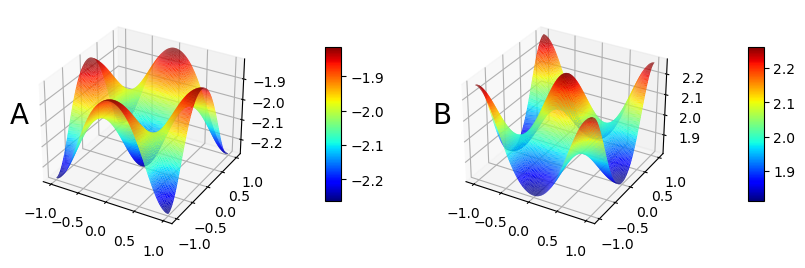}
  \includegraphics[width=2.2in]{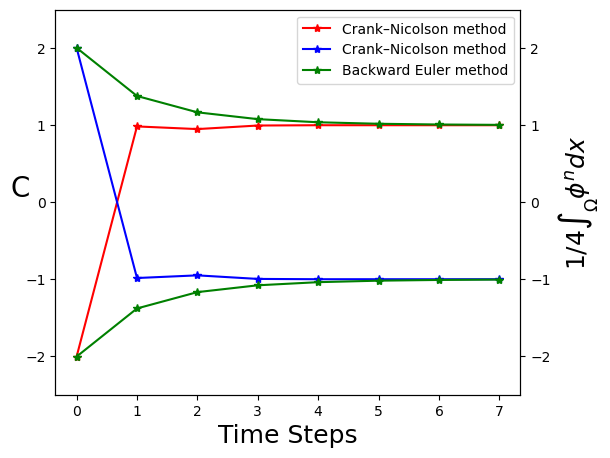}
  \caption{The solutions of $\phi^{n} \approx r + \delta B\cos(k\pi x_1)\cos(l\pi x_2)$ of the CN scheme for 2D function $\phi^{n+1}=c+\delta \cos(k\pi x_1)\cos(l\pi x_2)$ with $|\delta|=0.5$ are depicted in A and B panels. Here the perturbation function is $k=1, l=1$. The parameters are chosen as in A $(c=0.984375)$, ($r=-1.99310$) and B $(c=-0.984375)$, ($r=1.99310$), $\epsilon=0.1$, and $\Delta t=0.01$. This $r$ values are chosen from Table~\ref{Table_CN} in interval $(r_{1},r_{2})$. In panel C, the CN scheme jumps to a different solution with the initial conditions of $\phi^n$, ultimately converging to an incorrect solution. Conversely, the backward Euler scheme converges to a solution near $\phi^n$.}\label{compare1}
\end{figure}

Using $r + \psi(x)$ as an initial guess, we employ Newton's method to solve Eq.~\eqref{cn scheme} for $\phi^n$ given $\phi^{n+1}=c+\delta f(x)$. As an illustrative example, we present the results in Fig.~\ref{compare} {for the 1D case and Fig.~\ref{compare1} for the 2D case} with the parameters $c=0.984375$, $\epsilon=0.1$, and $\Delta t=0.01$. We initiate the process with $\delta=0.001$ and employ a homotopy continuation method to compute the solution with $\delta=0.5$ \cite{hao2022adaptive,hao2020homotopy,hao2020adaptive}. The solutions of $\phi^n$ corresponding to $\delta=0.5$ for different perturbation modes $k=1$ and $k=5$ are shown in Fig.~\ref{compare}. {For the 2D case, perturbation is given as $f(x)=\cos(\pi x_1)\cos(\pi x_2)$ and results are shown in Fig.~\ref{compare1}.
} In this particular instance, we choose $r=\pm 1.99310$ from the interval $\pm [r_1,r_2]$. 

Consequently, the CN scheme converges after a single time step but yields an incorrect solution. To elaborate, if we choose $\phi^n(x)\approx 1.99310$ as an initial condition, the CN scheme jumps to approximately $-1$ after a one-time step and continues to converge toward $-1$ after a few iterations. Conversely, the backward Euler scheme converges to a correct solution near $\phi^n$.{ Due to the robustness and stability of the backward Euler scheme, the numerical solution computed using this scheme serves as the reference solution for comparing with solutions computed using other schemes. Moreover, based on the PDE theory, the time evolution solution consistently converges to the nearby steady-state solution.  Thus we conclude that the CN scheme converges to an incorrect steady-state solution with this initial condition.}

\begin{figure}[!ht]
\centering
  \includegraphics[width=5in]{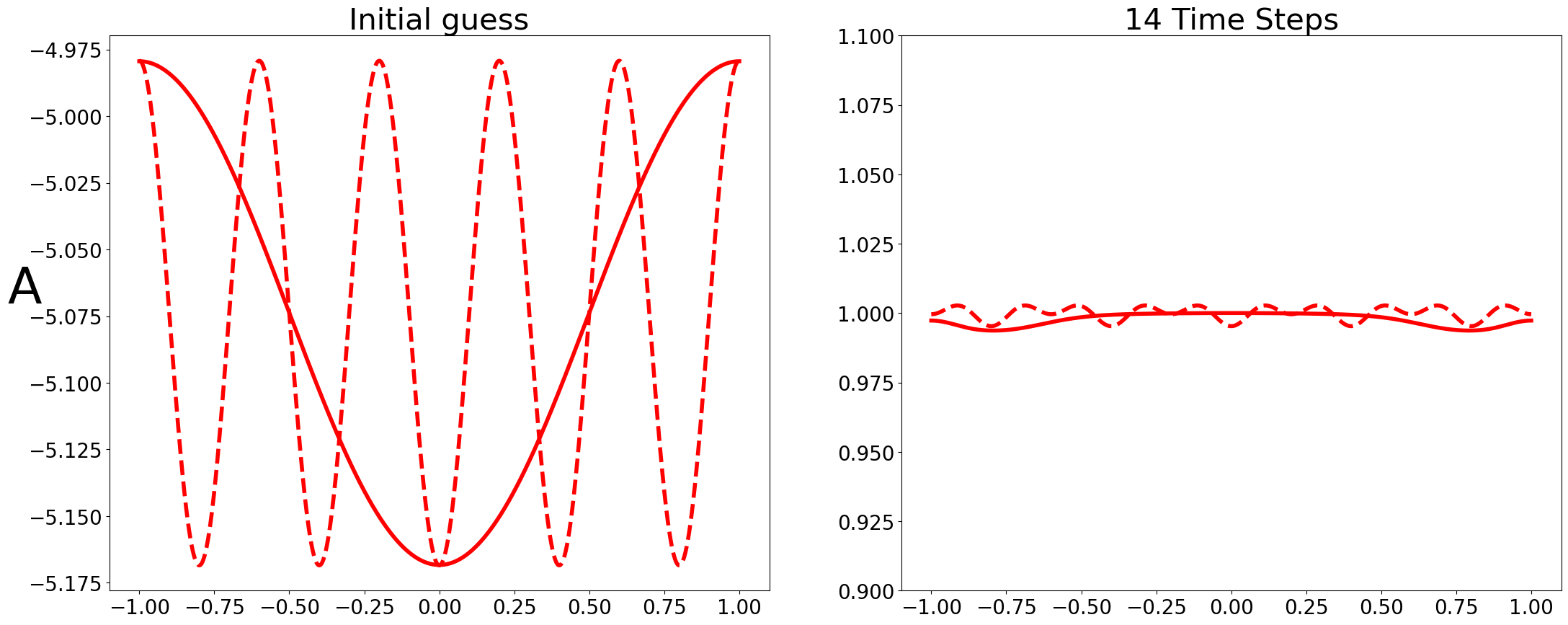}
  \includegraphics[width=5in]{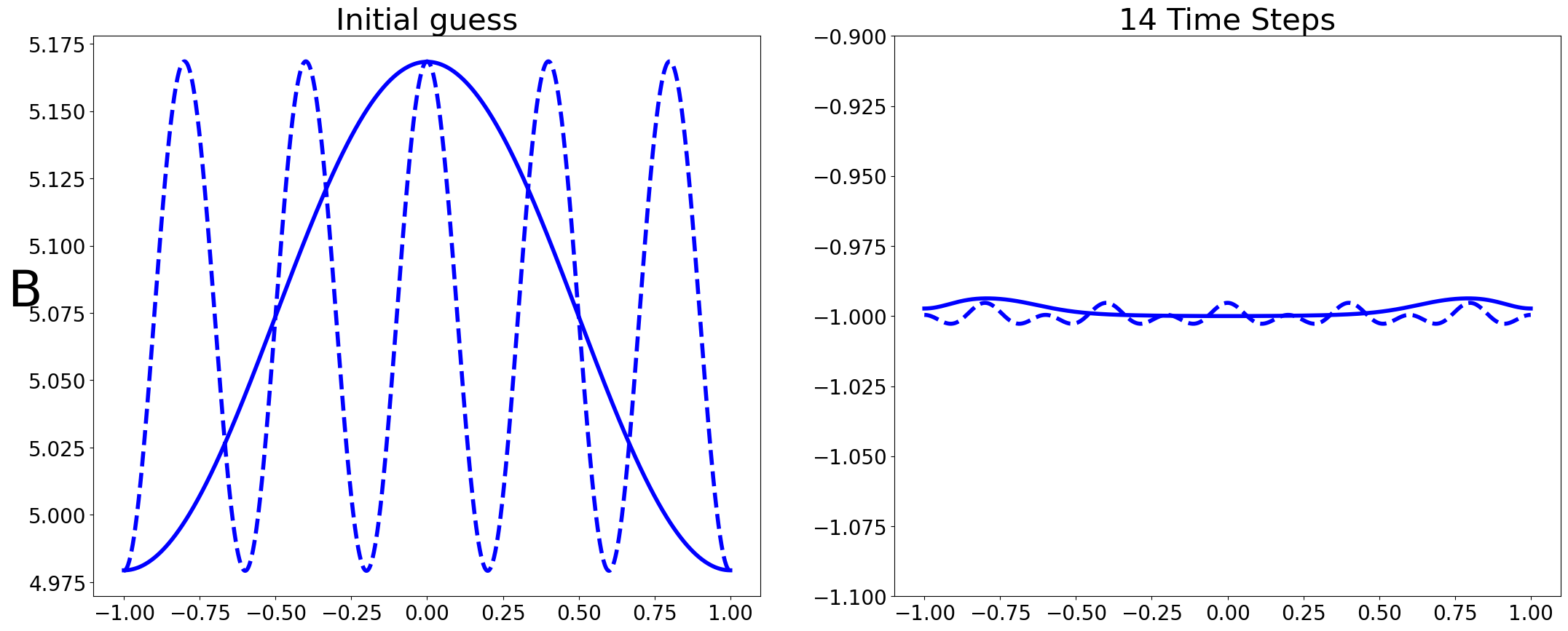}
    \includegraphics[width=4in]{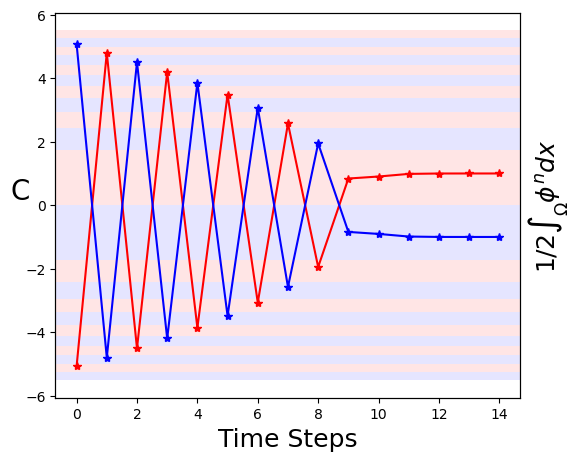}
  \caption{The solutions of $\phi^{n} \approx r + \delta B\cos(k\pi x_1)$ of the CN scheme for $\phi^{n+1}=c+\delta \cos(k\pi x_1)$ with $|\delta|=0.1$ are depicted in A and B panels. Here the solid curve is for $k=1$, the dashed curve is for $k=5$. The parameters are chosen as in A $(c=4.8)$, ($r=-5.074$) and B $(c=-4.8)$, ($r=5.074$), $\epsilon=0.1$, and $\Delta t=0.01$.  In panel C, the CN scheme jumps back and forth, ultimately converging to an incorrect solution. Red and blue regions are visualizing convergence intervals as Fig.~\ref{fig:rc}.
 } \label{PDEre}
\end{figure}

\begin{figure}[!ht]
\centering
  \includegraphics[width=5in]{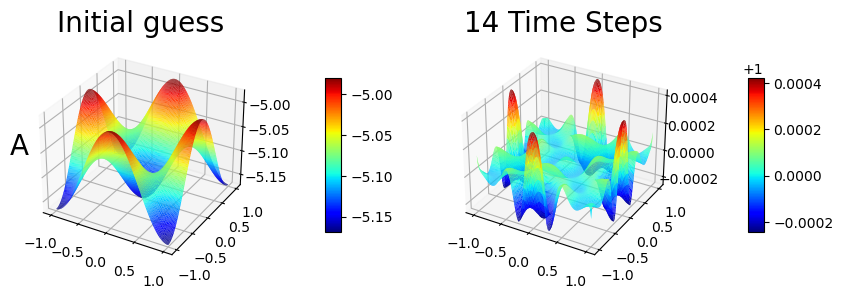}
  \includegraphics[width=5in]{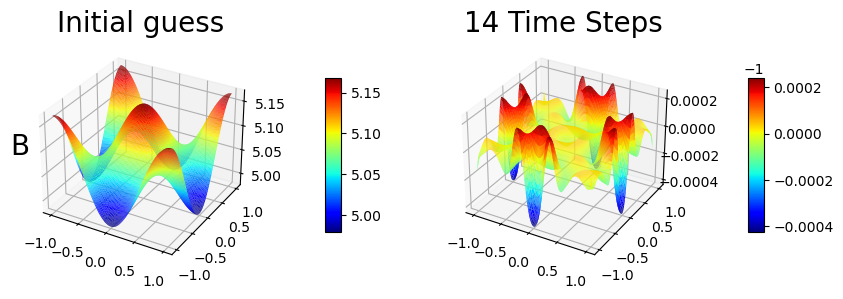}
    \includegraphics[width=4in]{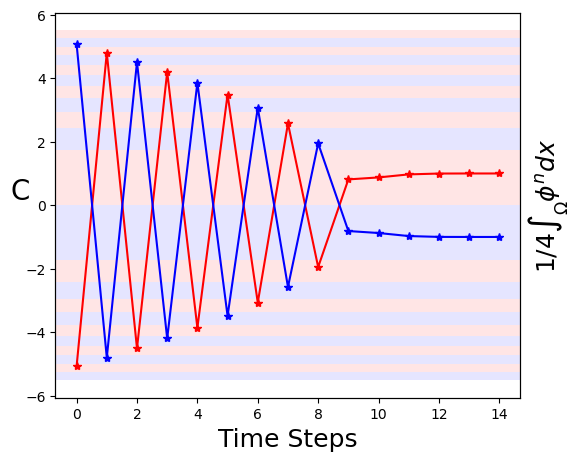}
  \caption{The solutions of $\phi^{n} \approx r + \delta B\cos(k\pi x_1)\cos(l\pi x_2)$ of the CN scheme for 2D function $\phi^{n+1}=c+\delta \cos(k\pi x_1)\cos(l\pi x_2)$ with $|\delta|=0.1$ are depicted in A and B panels. Here the perturbation function is $k=1, l=1$. The parameters are chosen as in A $(c=4.8)$, ($r=-5.074$) and B $(c=-4.8)$, ($r=5.074$) $\epsilon=0.1$, and $\Delta t=0.01$.  In panel C, the CN scheme jumps back and forth, ultimately converging to an incorrect solution. Red and blue regions are visualizing convergence intervals as Fig.~\ref{fig:rc}.
 } \label{PDEre1}
\end{figure}

Furthermore, we choose $r=\pm 5.074$ within the interval $\pm [r_9,r_{10}]$. As a result, the CN scheme converges after 9 zig-zag iterations but leads to an incorrect solution after 14 time steps shown in Fig.~\ref{PDEre} {for the 1D case and  Fig.~\ref{PDEre1} for the 2D case}. The zig-zag curves in Figs.~\ref{PDEre} and  \ref{PDEre1} demonstrate the patterns illustrated in Fig.~\ref{fig:rc}.
We choose different initial conditions with both $k=1$ and $k=5$ modes as well as $\delta=0.1$ for the 1D case in Fig.~\ref{PDEre} and for the 2D case in Fig.~\ref{PDEre1}.

\section{Convex Splitting of Modified Crank-Nicolson Scheme}
Next, we consider the convex splitting of the modified Crank-Nicolson (Mod CN) scheme of Eq.~\eqref{AC} defined as follows:

\begin{equation} \label{mod cn scheme}
\begin{cases}
        \frac{\phi^{n+1} - \phi^n}{\Delta t} -\frac{1}{2}(\triangle\phi^{n+1} + \triangle\phi^{n}) + \frac{1}{4\epsilon^2}\big( \phi^{n+1} + \phi^{n}\big) \big( (\phi^{n+1})^2 + (\phi^{n})^2\big)- \frac{1}{\epsilon^2}\big( \phi^{n}\big) = 0
    \\
    \frac{\partial \phi^n(x)}{\partial \mathbf{n}}= \frac{\partial \phi^{n+1}(x)}{\partial \mathbf{n}}=0~\hbox{~,~~~~~~} x\in \partial\Omega~.
\end{cases}
\end{equation}

\subsection{Unconditional stability}
We now examine the uniqueness of $\phi^{n+1} = c + \delta \psi(x)$ for a fixed $\phi^{n} = r + \delta f(x)$ within Eq.~\eqref{mod cn scheme}. Substituting these functions into Eq.~\eqref{mod cn scheme}, we obtain:
\begin{equation} 
   -\frac{1}{2}\triangle\psi+ \left[\frac{1}{\Delta t}+\frac{2c^2+(c+r)^2)}{4\epsilon^2}\right]\psi  +O(\delta)= G(x),\label{Eq_mcn}
\end{equation}
where $G(x)$ is given by:
\begin{equation} 
G(x)=\frac{   f(x)}{\Delta t}+\frac{1}{2} \triangle f(x)-\frac{1}{4\epsilon^2}\big( 3r^2f +c^2  f+2crf-4f\big)~.
\end{equation}
Given that $\frac{1}{\Delta t}+\frac{2c^2+(c+r)^2}{4\epsilon^2} > 0$, it becomes evident that only the trivial zero general solution exists within
\begin{equation} 
-\frac{1}{2}\triangle\psi+ \left[\frac{1}{\Delta t}+\frac{2c^2+(c+r)^2}{4\epsilon^2}\right]\psi = 0~.
\end{equation}

We conclude that the particular solution for Eq.~\eqref{Eq_mcn} is unique for any given $\phi^n$, independent of $\Delta t$. This demonstrates the unconditional stability of the convex splitting of the modified Crank-Nicolson scheme.

\subsection{The robustness analysis}
Similar to the CN scheme, we analyze the trivial solution structure of the convex splitting of the Mod CN  scheme firstly. Namely, for $\phi^{n+1}\equiv c$ and $\phi^{n}\equiv r$, we have
\begin{equation}
\frac{c - r}{\Delta t} + \frac{1}{4\epsilon^2}\big( r + c\big) \big( r^2 + c^2\big)- \frac{r}{\epsilon^2} = 0~.\label{mrc}
\end{equation}

Through a computation similar to the one used in analyzing the CN scheme \S \ref{3.2.1}, we can calculate the sequence of $r_n$ by solving Eq.~\eqref{mrc} with $c=r_{n-1}$ (here $r_0 = 0$). We present the values of $r_n$ for various ratios of $\frac{\Delta t}{2\epsilon^2}$  in Table~\ref{Table_mCN}.

\begin{table}[h]
\centering
\begin{tabular}{|c||c|c|c|c|}
\hline
$\frac{\Delta t}{2\epsilon^2}$ & $r_1$ & $r_2$ & $r_3$ & $r_4$\\
\hline
\hline
0.001 & 44.766 & 75.889 & 98.476 & 116.931 \\
\hline
0.01 & 14.283 & 24.165 & 31.334 & 37.192 \\
\hline
0.1 & 4.899 & 8.147 & 10.497 & 12.418 \\
\hline
0.25 & 3.464 & 5.641 & 7.212 & 8.497 \\
\hline
0.5 & 2.828 & 4.503 & 5.707 & 6.694 \\
\hline
\end{tabular}
\caption{Convergence Interval points $r_i$ for the Convex splitting of modified CN scheme with different $\frac{\Delta t}{2\epsilon^2}$.}\label{Table_mCN}
\end{table}

We can perturb the trivial solutions using $\phi^{n+1} = c + \delta f(x)$ and $\phi^{n} = r + \delta \psi(x)$ in a way similar to the one in \S \ref{PertAna}. This allows us to reformulate the Convex splitting of Mod CN scheme as follows:
\begin{equation} \label{mod parti_equ}
    -\frac{1}{2} \triangle\psi+\left(-\frac{1}{\Delta t}+\frac{3r^2+c^2+2cr-4}{4\epsilon^2}\right)\psi  = G(x)~,
\end{equation}
where $G(x)$ is given by:
\begin{equation} 
G(x)=-\frac{   f(x)}{\Delta t}+\frac{1}{2} \triangle f(x)-\frac{1}{4\epsilon^2}\big( 3c^2f +r^2  f+2crf\big)~.
\end{equation}
When we choose $f(x) = \cos(k\pi x_1)$, the particular solution becomes $\psi = B \cos(k \pi x_1)$ {for the 1D case}, where the coefficient $B$ is determined by the following expression from using Eq.~\eqref{mod parti_equ}:
\begin{equation}
    B = - \frac{\frac{1}{\Delta t}+\frac{k^2 \pi^2}{2}+\frac{3c^2+r^2+2cr}{4 \epsilon^2}}{-\frac{1}{\Delta t}+\frac{3r^2+c^2+2cr-4}{4
\epsilon^2}+\frac{k^2 \pi^2}{2}}=- \frac{\frac{3c^2+r^2+2cr}{2 \epsilon^2}+\frac{2}{\Delta t}+{k^2 \pi^2}}{\frac{3r^2+c^2+2cr-4}{2
\epsilon^2}-\frac{2}{\Delta t}+{k^2 \pi^2}}.
\end{equation}

Using the computed $\psi(x)$ as an initial approximation, we employ Newton's method to solve for $\phi^n$ given $\phi^{n+1}=c+\delta f(x)$ in Eq.~\eqref{mod cn scheme}. This results in a solution structure similar to that of the CN scheme.

For 2D we have similar results as $f(x) = \cos(k\pi x_1) \cos(l \pi x_2)$ for $k, l \in \mathbb{N}^+$. We find that $\psi = B \cos(k \pi x_1) \cos(l \pi x_2)$, with the coefficient $B$ given by 
\begin{equation}
    B = - \frac{\frac{3c^2+r^2+2cr}{2 \epsilon^2}+\frac{2}{\Delta t}+{k^2 \pi^2}+{l^2 \pi^2}}{\frac{3r^2+c^2+2cr-4}{2
\epsilon^2}-\frac{2}{\Delta t}+{k^2 \pi^2}+{l^2 \pi^2}}.
\end{equation}

\section{Diagonally Implicit Runge–Kutta (DIRK)}
The DIRK family of methods is the most widely used implicit Runge-Kutta (IRK) method for solving phase field modeling problems due to their relative ease of implementation \cite{kennedy2016diagonally}.
{Some applications of the DIRK methods on the Allen-Cahn equation can be found in \cite{church2019high, shah2018efficient}}. These methods are characterized by a lower triangular A-matrix with at least one non-zero diagonal entry and are sometimes referred to as semi-implicit or semi-explicit Runge-Kutta methods. This structure allows for solving each stage individually rather than all stages simultaneously. We write the general formula of DIRK method in Butcher array format as follows:
\begin{center}
\begin{tabular}{ c| c c c c c c } 

 $c_1$ & $a_{11}$ & 0&0&0&\dots&0 \\ 
 $c_2$ & $a_{21}$ & $a_{22}$&0&0&\dots&0 \\ 
 $c_3$ & $a_{31}$ & $a_{32}$&$a_{33}$&0&\dots&0 \\ 
 $c_3$ & $a_{41}$ & $a_{42}$&$a_{43}$&$a_{44}$&\dots&0 \\
  $\vdots$ &$\vdots$ & $\vdots$&$\vdots$&$\vdots$&$\ddots$&$\vdots$ \\
 $c_s$ & $a_{s1}$ & $a_{s2}$&$a_{s3}$&$a_{s4}$&\dots&$a_{ss}$ \\
 \hline
   & $b_{1}$ & $b_{2}$&$b_{3}$&$b_{4}$&\dots&$b_{s}$ \\ 
\end{tabular}.
\end{center}
When solving the Allen-Cahn equation, the DIRK method is summarized as 
\begin{equation}
    \phi^{n+1}=\phi^n+\triangle t \sum_{i=1}^s b_i k_i,
\end{equation}
where
\begin{equation}
    k_i=F(\phi^n + \triangle t \sum_{j=1}^i a_{ij}k_j)\hbox{~and~}    F(\phi)=\triangle \phi -\frac{1}{\epsilon^2}(\phi^3-\phi).
\end{equation}
By letting $\phi^{n+1}_i=\phi^n + \triangle t \sum_{j=1}^i a_{ij}k_j$, we can rewrite the DIRK method as \cite{kennedy2016diagonally, kennedy2019diagonally}
\begin{equation}
\phi^{n+1}_i =
\begin{cases}
    \phi^n & \text{for } i = 0, \\
    \phi^{n+1}_0 + \Delta t \sum_{j=1}^i a_{ij} F(\phi^{n+1}_j) & \text{for } 1 \leq i \leq s.
\end{cases}
\end{equation}
In this case, the final solution can be expressed as $    \phi^{n+1}=\phi^{n+1}_0 + \Delta t \sum_{i=1}^s b_i F(\phi^{n+1}_i)$. Then the DIRK method represents a multi-stage backward Euler method which solves for $\phi^{n+1}_i$ for each stage $i$.

\subsection{The stability condition via bifurcation analysis}\label{scvba}
First, we consider a trivial solution case of  $\phi^{n+1}_{i+1} \equiv c$ for given $\phi^{n+1}_0, \phi^{n+1}_1, \dots, \phi^{n+1}_{i-1}$, specifically, $\phi^{n+1}_i  \equiv r$:
\begin{equation}\label{fis2st}
    \phi^{n+1}_{i+1}=\phi^{n+1}_0 + \Delta t \sum_{j=1}^{i+1} a_{i+1,j} F(\phi^{n+1}_j)
\end{equation}
Next, let's perturb  the trivial solutions with $\phi^{n+1}_{i+1}(x) = c + \delta \psi(x)$ and $\phi^{n+1}_{i}(x) = r + \delta f(x)$. By substituting into Eq.~\eqref{fis2st} and retaining the linear term in $\delta$, we obtain

\begin{equation}
-\triangle \psi(x)+\Big(\frac{1}{\triangle t  a_{i+1,i+1}}+\frac{3c^2-1}{\epsilon^2}\Big)\psi(x)  =
\begin{cases}
   \frac{f(x)}{a_{1,1}}  +O(\delta) & $ if $ i=0\\
      \frac{a_{i+1,i}}{a_{i+1,i+1}}\Big( \triangle f-\frac{1}{\epsilon^2}(3 r^2 f - f )\Big) +O(\delta) & $ else$
\end{cases}    
\end{equation}

We aim to assess the uniqueness of $\phi^{n+1}_{i+1}$ while holding $\phi^{n+1}_{i}$ fixed. In this case, $f(x)$ is a given function, and we focus on the uniqueness of the homogeneous solution, namely,
\begin{equation}
-\triangle \psi(x)+\Big(\frac{1}{\triangle t  a_{i+1,i+1}}+\frac{3c^2-1}{\epsilon^2}\Big)\psi(x) = 0.
\end{equation}
This is essentially the same as Eq.~\eqref{Linear1} if $a_{i+1,i+1}=1$. By {\bf Proposition \ref{prop1}}, we have that the solution is unique if
\begin{equation}
    \Big(\frac{1}{\triangle t  a_{i+1i+1}}+\frac{3c^2-1}{\epsilon^2}\Big)\ge 0,
\end{equation}
or
\begin{equation}
    0\le \triangle t a_{i+1i+1}\le \epsilon^2.
\end{equation}
Consequently, we can conclude that the particular solution for Eq.~\eqref{fis2st} is unique, and $\phi^{n+1}_{i+1}$ is also unique while satisfying the stability condition $\Delta t a_{i+1i+1} \le \epsilon^2$. Then the stability condition of the DIRK method is
\begin{equation}
    0 \le \triangle t \le\frac{\epsilon^2}{\max_i a_{ii}}.
\end{equation}

\subsection{The robustness analysis}

In theory, we can apply robustness analysis to any order of the DIRK method. In this section, for simplicity, we illustrate the idea by considering the 2nd order DIRK method with the following \(2\times2\)  Butcher array:
\[
\begin{array}{ c| c c } 
1/4 & 1/4 & 0\\ 
3/4 & 1/2 & 1/4\\ 
\hline
 & 1/2 & 1/2\\ 
\end{array}
\]
We first analyze the trivial solution case, namely \(\phi^{n+1}\equiv c\) and \(\phi^{n}\equiv r\).
 Then, the DIRK method with 2nd order for solving the Allen-Cahn equation is expressed as:
\begin{equation}
\begin{cases}
    \phi^{n+1}_1= \phi^{n+1}_0+\frac{\triangle\phi^{n+1}_1}{4}-\frac{\Delta t}{4\epsilon^2}((\phi^{n+1}_1)^3-\phi^{n+1}_1),\\
    \phi^{n+1}_2= \phi^{n+1}_1+\frac{\triangle\phi^{n+1}_1}{4}-\frac{\Delta t}{4\epsilon^2}((\phi^{n+1}_1)^3-\phi^{n+1}_1)+\frac{\triangle\phi^{n+1}_2}{4}-\frac{\Delta t}{4\epsilon^2}((\phi^{n+1}_2)^3-\phi^{n+1}_2),\\
    \phi^{n+1}= \phi^{n+1}_2+\frac{\triangle\phi^{n+1}_2}{4}-\frac{\Delta t}{4\epsilon^2}((\phi^{n+1}_2)^3-\phi^{n+1}_2). \label{DIRK2}
\end{cases}    
\end{equation}
By letting $\phi^{n+1}\equiv 0$, we have 3 roots for $\phi^{n+1}_2$ as $\pm \frac{r_1}{2}$
 and $ 0$ by solving the last equation in Eq.~\eqref{DIRK2}, where $r_1=2\sqrt{{1+\frac{4 \epsilon^2}{\triangle t}}}$. By simplifying Eq.~\eqref{DIRK2} with $\phi^{n+1}\equiv 0$, we have 
\begin{equation}\label{finpoly} 2\phi^{n+1}_1-2\phi^{n+1}_2=\phi^{n+1}_0.
\end{equation}
By plugging Eq.~\eqref{finpoly} into the first equation of Eq.~\eqref{DIRK2}, we have 
\begin{equation}\label{finpoly1} 
\frac{-\phi^{n+1}_0+2\phi^{n+1}_2}{2}=-\frac{\Delta t}{4\epsilon^2}\left(\left(\frac{\phi^{n+1}_0+2\phi^{n+1}_2}{2}\right)^3-\frac{\phi^{n+1}_0+2\phi^{n+1}_2}{2}\right).
\end{equation}

If $\phi^{n+1}_2=\pm \frac{r_1}{2}$, since of the discriminant of the cubic polynomial from the second equation in Eq.~\eqref{DIRK2},
\begin{equation}
-27\left(\frac{8\epsilon^2}{\Delta t}\phi^{n+1}_2\right)^2+4\left(\frac{4\epsilon^2}{\Delta t}+1\right)^3<0, \quad \text{ when } \Delta t < 4 \epsilon^2,
\end{equation}
we have only one root as $\phi^{n+1}_1=\pm h_1$. Thus we have
\begin{equation}
    \phi^{n+1}_2=\pm \frac{r_1}{2}, \quad \phi^{n+1}_1=\pm \frac{-s_1+2\phi^{n+1}_2}{2},\hbox{~and~} \phi^{n+1}_0=\mp s_1,
\end{equation}
where $s_1$ is the root of
\begin{equation}
\frac{s_1+r_1}{2}=-\frac{\Delta t}{4\epsilon^2}\left(\left(\frac{-s_1+r_1}{2}\right)^3-\frac{-s_1+r_1}{2}\right).
\end{equation}

If $\phi^{n+1}_2=0$, we have another set of solutions:
\begin{equation}
    \phi^{n+1}_2=0, \quad \phi^{n+1}_1=\pm \frac{r_1}{2}\text{ or } 0, \quad \phi^{n+1}_0=2\phi^{n+1}_1-2\phi^{n+1}_2.
\end{equation}

For $ 0< r_1<s_1$, we have five roots for $\phi^{n+1}_0$, namely, $\phi^{n+1}_0=\pm r_1,\mp s_1$, and $0$.

\begin{figure}[h]
    \centering
  \includegraphics[width=6in]{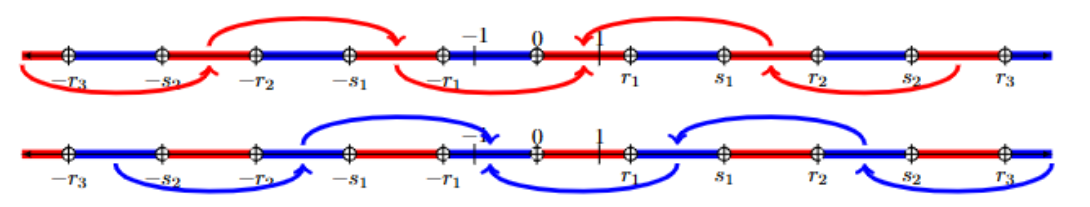}
    \caption{Visualizing Convergence Intervals of DIRK scheme with 2nd order in Eq.~\eqref{DIRK2}. Suppose initial conditions are chosen in red regions. In that case, DIRK with 2nd order scheme eventually converges to $1$, while the initial conditions chosen in blue regions lead the DIRK method to converge to $-1$. The values of $r_i, s_i$ with different $\frac{\Delta t}{4 \epsilon^2}$ are shown in Table~\ref{Table_rk}.}
    \label{fig:rkc}
\end{figure}

By letting {$\phi^{n+1} =  s_1$}, we obtain a unique solution for $\phi^n$ due to the discriminant of the cubic polynomial, denoted as $\phi^n = s_2$. Inductively, we can define $s_i$ and $r_i$ by solving Eq.~\eqref{DIRK2} for $\phi^n$ with $\phi^{n+1} = s_{i-1}$ and $\phi^{n+1} = r_{i-1}$, respectively. The values of $r_i$ and $s_i$ for different values of $\frac{\Delta t}{4\epsilon^2}$ are shown in Table~\ref{Table_rk}, and the iterations of the DIRK method in different regions are illustrated in Fig.~\ref{fig:rkc}.

\begin{table}[h]
\centering
\begin{tabular}{|c||c|c|c|c|c|c|c|c|}
\hline
$\frac{\Delta t}{4\epsilon^2}$ & $r_1$ & $s_1$ & $r_2$ & $s_2$& $r_3$ & $s_3$ & $r_4$ & $s_4$\\
\hline
\hline
0.001& 63.277& 159.524& 280.251& 421.311& 580.137& 754.936& 944.371& 1147.391\\
0.01& 20.1& 50.612& 88.857& 133.527& 183.81& 239.141& 299.098& 363.349\\
0.1& 6.633& 16.517& 28.821& 43.14& 59.221& 76.889& 96.012& 116.485\\
0.25& 4.472& 10.958& 18.95& 28.2& 38.552& 49.898& 62.156& 75.262\\
0.5& 3.464& 8.306& 14.188& 20.942& 28.462& 36.675& 45.524& 54.966 \\
\hline
\end{tabular}
\caption{Convergence Interval points $r_i$, $s_i$ for the DIRK 2nd order scheme in Eq.~\eqref{DIRK2} with different $\frac{\Delta t}{4\epsilon^2}$.}\label{Table_rk}
\end{table}

Next, we perturb the trivial solutions using $\phi^{n+1} = c + \delta \xi (x)$, $\phi^{n+1}_2=c_2+\delta \xi_2(x)$, $\phi^{n+1}_1=c_1+\delta \xi_1(x)$ and $\phi^{n+1}_0 = r + \delta \xi_{0}(x)$ as similar way in \S \ref{PertAna}, where $\xi(x)$ is a given perturbed function. After plugging in to Eq.~\eqref{DIRK2} and retaining the linear term in $\delta$, we have
\begin{equation}
\begin{cases}
     \xi_1= \xi_0+\frac{\triangle t}{4} \triangle \xi_1 -\frac{\triangle t}{4\epsilon^2}(3 \xi_1(c_1)^2- \xi_1)+O(\delta),\\
    \xi_2= \xi_1+\frac{\triangle t}{4} \triangle \xi_2 +\frac{\triangle t}{4} \triangle \xi_1 -\frac{\triangle t}{4\epsilon^2}(3 \xi_1(c_1)^2- \xi_1)-\frac{\triangle t}{4\epsilon^2}(3 \xi_2(c_2)^2- \xi_2)+O(\delta),\\
    \xi (x)= \xi_2+\frac{\triangle t}{4} \triangle \xi_2-\frac{\triangle t}{4\epsilon^2}(3\xi_2(c_2)^2- \xi_2)+O(\delta).
\end{cases}
\end{equation}

By choosing specific functions {in the 1D case} as
\begin{equation}
    \xi(x)=\cos( k \pi x_1), \quad \xi_2 (x)=B_2\cos(k \pi  x_1), \quad \xi_1 (x)=B_1\cos(k \pi  x_1), \quad \xi_0 (x)=B_0\cos(k \pi  x_1)
\end{equation}
we obtain
\begin{equation}
\begin{cases}
    B_2=\frac{1}{1-\frac{\Delta t k^2\pi^2}{4}-\frac{\Delta t}{4 \epsilon^2}(3(c_2)^2-1)}\\
    B_1= B_2\frac{1+\frac{\Delta t k^2\pi^2}{4}+\frac{\Delta t}{4 \epsilon^2}(3(c_2)^2-1)}{1-\frac{\Delta t k^2\pi^2}{4}-\frac{\Delta t}{4 \epsilon^2}(3(c_1)^2-1)}\\
    B_0= B_1(1+\frac{\Delta t k^2\pi^2}{4}+\frac{\Delta t}{4 \epsilon^2}(3(c_1)^2-1))
\end{cases}
\end{equation}

\begin{figure}[h]
\centering
  \includegraphics[width=2in]{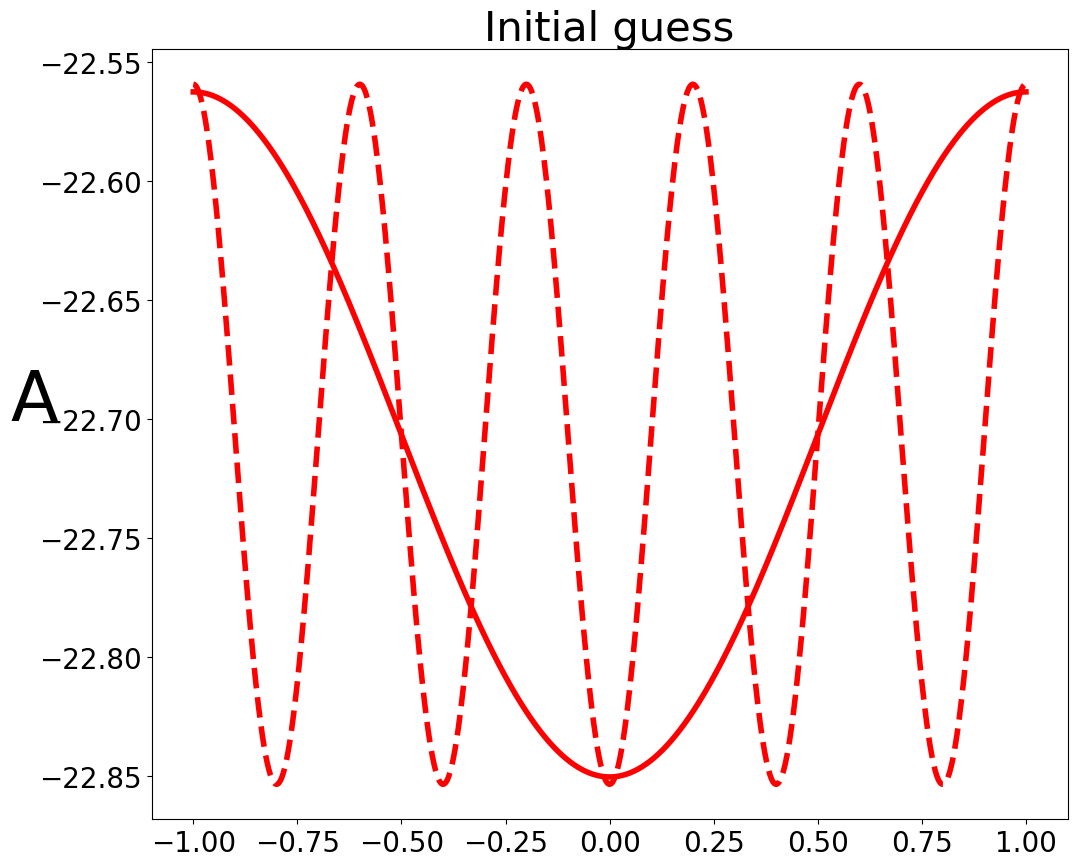}
  \includegraphics[width=2in]{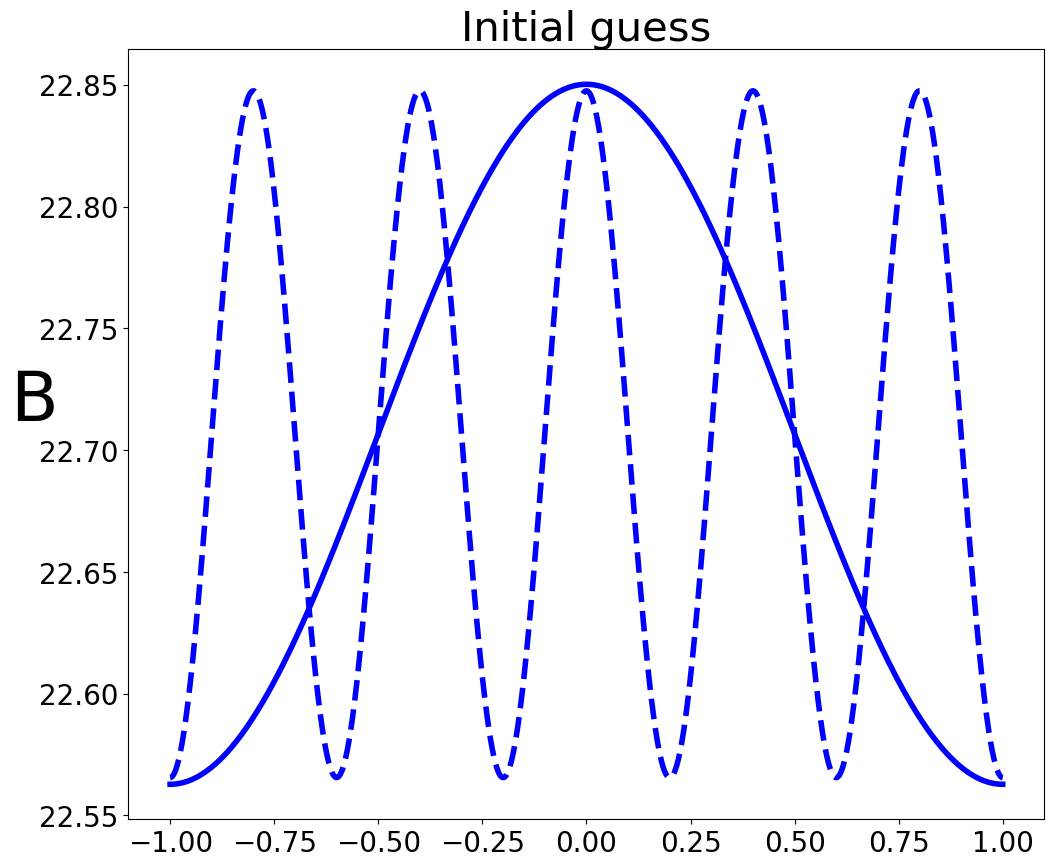}
  \includegraphics[width=2.2in]{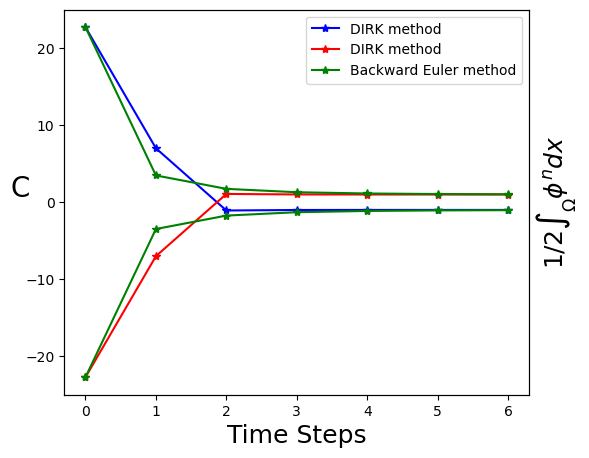}
  \caption{The solutions of $\phi^{n} \approx r + \delta B_0\cos(k\pi x_1)$ of the DIRK scheme for $\phi^{n+1}=c+\delta \cos(k\pi x_1)$ with $|\delta|=0.1$ are depicted in A and B panels. Here the solid curve is for $k=1$, the dashed curve is for $k=5$. The parameters are chosen as in A $(c=-7)$, ($r=-22.70665$) and B $(c=7)$, ($r=22.70665$), $\epsilon=0.1$, and $\Delta t=0.01$. This $c$ values are chosen from Table~\ref{Table_rk} in interval $(r_{1},s_{1})$. In panel C, the DIRK scheme jumps to a different solution with the initial conditions of $\phi^n$, ultimately converging to an incorrect solution. Conversely, the backward Euler scheme converges to a solution near $\phi^n$.}\label{PDErk11}
\end{figure}

\begin{figure}[h]
\centering
  \includegraphics[width=4in]{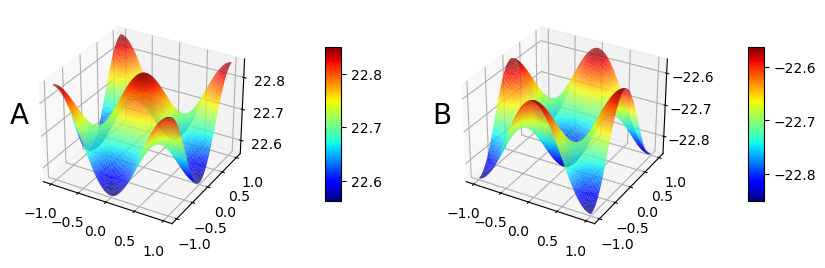}
  \includegraphics[width=1.5in]{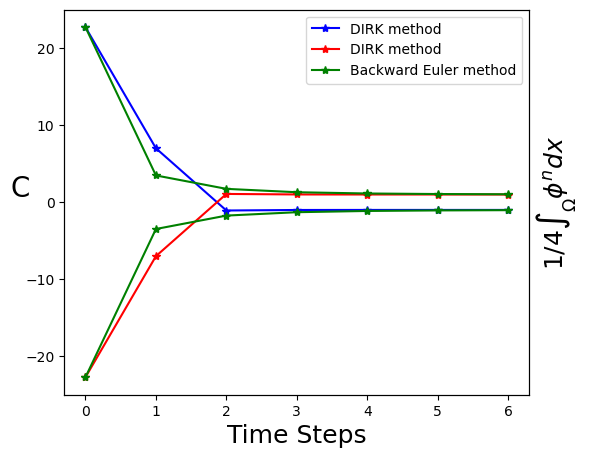}
  \caption{The solutions of $\phi^{n} \approx r + \delta B_0\cos(k\pi x_1)\cos(l\pi x_2)$ of the DIRK scheme for 2D function $\phi^{n+1}=c+\delta \cos(k\pi x_1)\cos(l\pi x_2)$ with $|\delta|=0.1$ are depicted in A and B panels. Here the perturbation function is $k=1, l=1$. The parameters are chosen as in A $(c=-7)$, ($r=-22.70665$) and B $(c=7)$, ($r=22.70665$), $\epsilon=0.1$, and $\Delta t=0.01$. This $c$ values are chosen from Table~\ref{Table_rk} in interval $(r_{1},s_{1})$. In panel C, the DIRK scheme jumps to a different solution with the initial conditions of $\phi^n$, ultimately converging to an incorrect solution. Conversely, the backward Euler scheme converges to a solution near $\phi^n$.}\label{PDErk112d}
\end{figure}

Then, we employ Newton's method to solve Eq.~\eqref{DIRK2} for $\phi^n$,  given $\phi^{n+1}=c+\delta \xi(x)$, taking the initial guess $\phi^{n} \approx r + \delta B_0\cos(k\pi x_1)$. As an illustrative example, we show the solutions of $\phi^n$ in Fig.~\ref{PDErk11} for the 1D case and Fig.~\ref{PDErk112d} for the 2D case with the parameters $c=-7$, $\epsilon=0.1$, and $\Delta t=0.01$. We initiate the process with $\delta=0.001$ and employ a homotopy continuation method to compute the solution with $\delta=0.1$ \cite{hao2022adaptive, hao2020homotopy,hao2020adaptive}. The solutions of $\phi^n$ corresponding to $\delta=0.1$ for different perturbation modes $k=1$ and $k=5$ are shown in Fig.~\ref{PDErk11}. {For the 2D case, perturbation is given as $\xi(x)=\cos(\pi x_1)\cos(\pi x_2)$ and we present the results in Fig.~\ref{PDErk112d}.}

Thus we use $\phi^n$ as the initial condition to solve the Allen-Cahn equation by using DIRK with second order. Since we choose $c=\pm 7$ from the interval $\pm [r_1,s_1]$, consequently, the DIRK scheme converges after two time steps but yields an incorrect solution. To elaborate, if we choose $\phi^{n+1}(x)\approx 7$ as an initial condition, the DIRK scheme jumps to approximately $-1$ after a one-time step and continues to converge toward $-1$ after a few iterations. Conversely, the backward Euler scheme converges to a correct solution near $\phi^n$. 

\begin{figure}[!th]
\centering
  \includegraphics[width=5in]{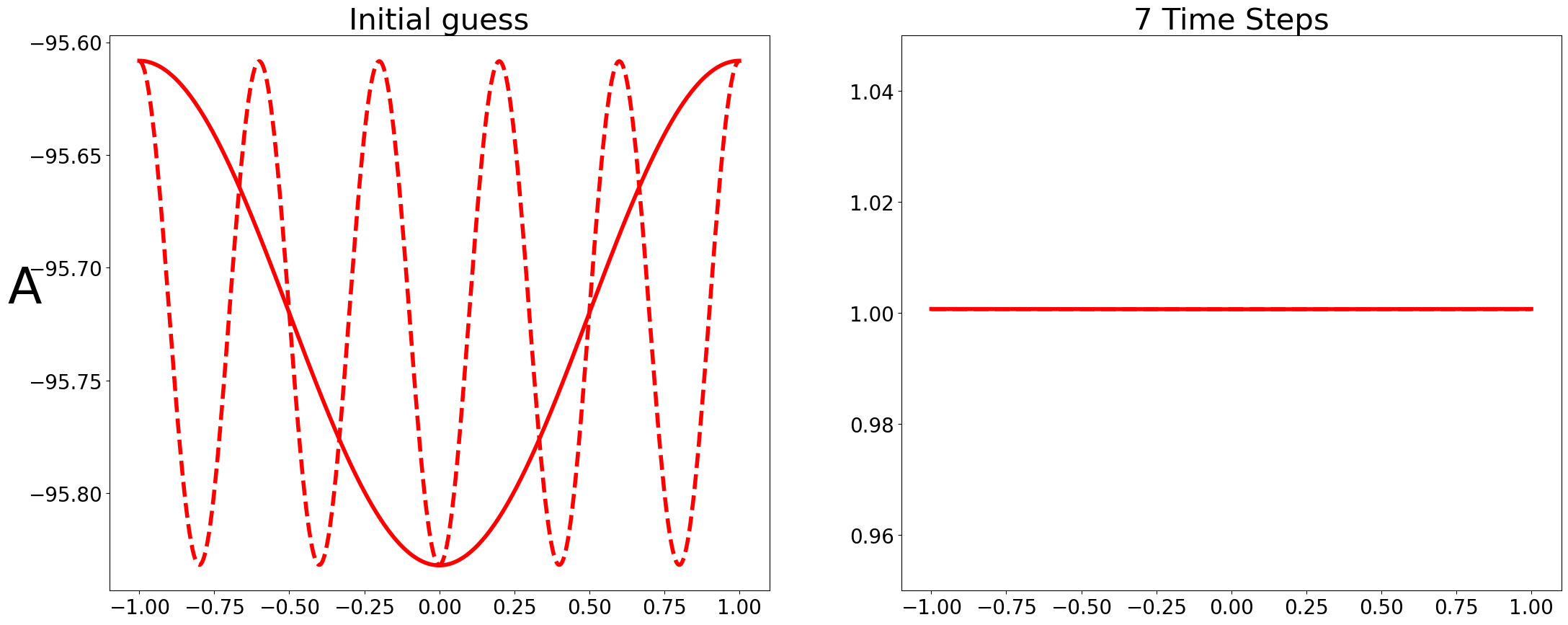}
  \includegraphics[width=5in]{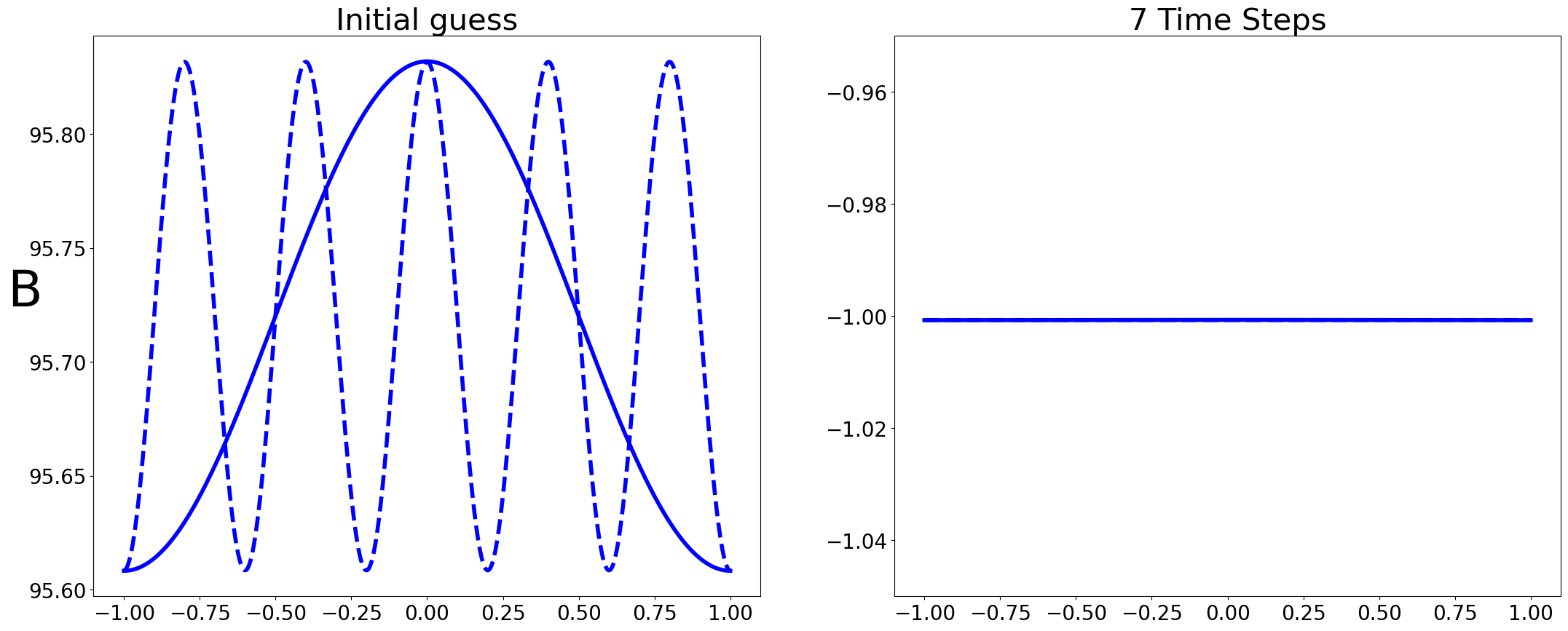}
    \includegraphics[width=4in]{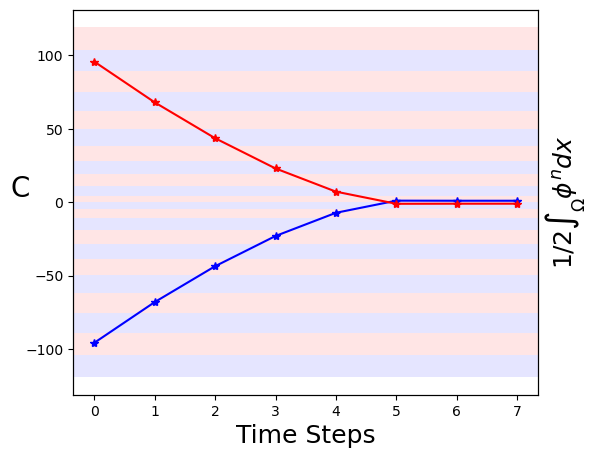}
  \caption{The solutions of $\phi^{n} \approx r + \delta B_0\cos(k\pi x_1)$ of the DIRK scheme for  function $\phi^{n+1}=c+\delta \cos(k\pi x_1)$ with $|\delta|=0.1$ are depicted in A and B panels. Here the solid curve is for $k=1$, the dashed curve is for $k=5$. The parameters are chosen as in A $(c=-68)$, ($r=-95.72$) and B $(c=68)$, ($r=95.72$), $\epsilon=0.1$, and $\Delta t=0.01$.  In panel C, the DIRK scheme, ultimately converging to an incorrect solution. Red and blue regions are visualizing convergence intervals as Fig.~\ref{fig:rkc}.
 } \label{PDErk1}
\end{figure}

\begin{figure}[!th]
\centering
  \includegraphics[width=5in]{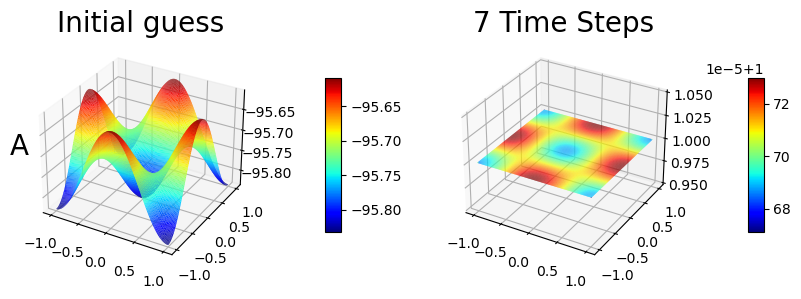}
  \includegraphics[width=5in]{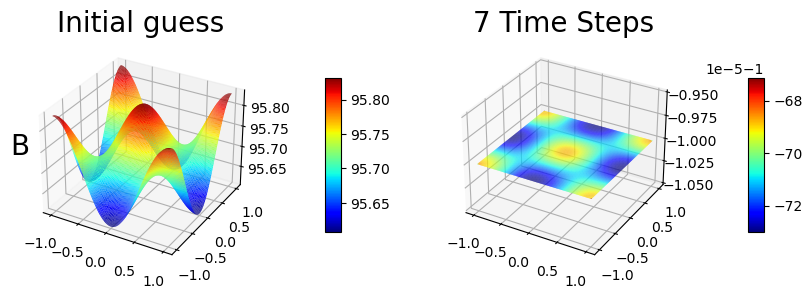}
    \includegraphics[width=4in]{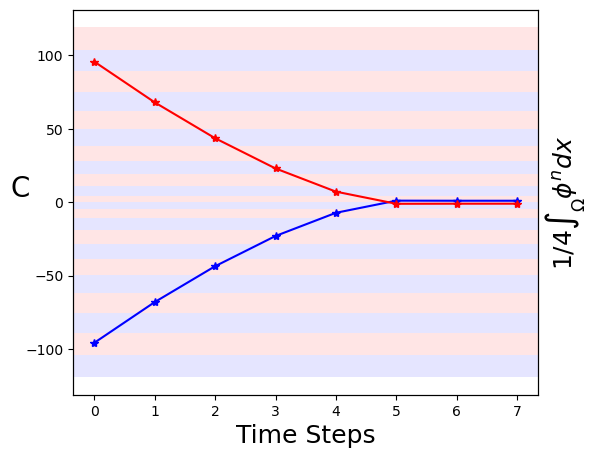}
  \caption{The solutions of $\phi^{n} \approx r + \delta B_0\cos(k\pi x_1)\cos(l\pi x_2)$ of the DIRK scheme for  2D function $\phi^{n+1}=c+\delta \cos(k\pi x_1)\cos(l\pi x_2)$ with $|\delta|=0.1$ are depicted in A and B panels. Here the perturbation function is $k=1, l=1$. The parameters are chosen as in A $(c=-68)$, ($r=-95.72$) and B $(c=68)$, ($r=95.72$), $\epsilon=0.1$, and $\Delta t=0.01$.  In panel C, the DIRK scheme, ultimately converging to an incorrect solution. Red and blue regions are visualizing convergence intervals as Fig.~\ref{fig:rkc}.
 } \label{PDErk2}
\end{figure}

Furthermore, we choose $r=\pm 95.72$ within the interval $\pm [r_5,s_{5}]$. As a result, the DIRK scheme converges after five iterations but leads to an incorrect solution. This convergence process is evident in the jumping curves depicted in Fig.~\ref{fig:rkc}. The initial conditions for $k=1$ and $k=5$ modes with $\delta=0.1$, as well as the final solutions after 7 time steps, are shown in Fig.~\ref{PDErk1} for the 1D case and Fig.~\ref{PDErk2} for the 2D case.


\section{Conclusions}
The Allen-Cahn equation serving as a fundamental tool for modeling phase transitions, offers invaluable insights into interface evolution across diverse physical systems. In this paper, we have devoted into the stability and robustness of various time-discretization numerical schemes utilized to solve the Allen-Cahn equation, recognizing their pivotal role in ensuring precise simulations in practical applications.

Our stability analyses of several numerical methods, including the backward Euler, Crank-Nicolson,  Convex Splitting of modified Crank-Nicolson schemes, and the DIRK method, have unveiled fundamental stability conditions for each method. Notably, the backward Euler scheme, Crank-Nicolson, and DIRK methods exhibited conditional stability, necessitating careful consideration of time step sizes. Conversely, the convex splitting of the modified Crank-Nicolson scheme showcased unconditional stability, affording flexibility in time step selection without compromising numerical accuracy. 

\begin{table}[h]
\centering
\begin{tabular}{|c||c|c|c|c|}
\hline
Numerical Scheme & Backward Euler & CN & Convex Splitting of Modified CN & DIRK\\
\hline
\hline
Stability Condition & $\Delta t \le \epsilon^2$ & $\Delta t \le 2\epsilon^2$ & $\Delta t \le \infty$ & $\Delta t \le \frac{\epsilon^2}{max_i a_{ii}}$ \\
\hline

\end{tabular}
\caption{Summary of stability conditions for different numerical schemes. $a_{ii}$ in DIRK is the diagonal elements in Butcher array format.}\label{Table_all}
\end{table}

Furthermore, our robustness analyses have shed light on the behavior of numerical solutions under varying initial conditions. While the backward Euler method demonstrated robustness, reliably converging to physical solutions regardless of initial conditions; other methods such as the Crank-Nicolson and convex splitting of modified Crank-Nicolson schemes, as well as the DIRK method, exhibited sensitivity to initial conditions in the solving of these nonlinear schemes at each step, potentially leading to wrong solutions if the initial conditions are not carefully chosen.

In conclusion, our study introduces the concepts of stability and robustness to the realm of numerical methods for solving the Allen-Cahn equation. By elucidating the stability conditions and robustness characteristics of these methods, we provide a novel framework for evaluating numerical techniques tailored to nonlinear differential equations, thereby advancing the accuracy and reliability of phase transition simulations in various scientific domains.

\section{Acknowledgement}
SL and WH are supported by NIH via 1R35GM146894.

\bibliographystyle{plain}
\bibliography{ref}

\end{document}